\documentclass[10pt,twocolumn,superscriptaddress,nofootinbib]{revtex4-1}
\usepackage[total={7.0in,9.5in}, top=1.0in, includefoot]{geometry}
\usepackage{epsfig}
\usepackage{subfigure}
\usepackage{amsmath}
\usepackage[usenames,dvipsnames]{color}
\usepackage{amssymb}
\usepackage{setspace}
\usepackage{graphicx} 
\usepackage{times}
\usepackage{float}
\usepackage{amsthm}
\usepackage{hyperref}
\usepackage{multirow}
\usepackage{array}
\usepackage{xspace}
\usepackage{cprotect}
\usepackage{diagbox}
\hypersetup{bookmarks=true, unicode=false, pdftoolbar=true, pdfmenubar=true, pdffitwindow=false, pdfstartview={FitH}, pdfcreator={Robert Steele}, pdfproducer={Daniel Larremore}, pdfkeywords={} {} {}, pdfnewwindow=true, colorlinks=true, linkcolor=red, citecolor=Green, filecolor=magenta, urlcolor=cyan,}
\usepackage{enumitem}

\makeatletter
\newcommand{\px}{\ext@arrow 0359\rightarrowfill@{}{X}}
\makeatother

\makeatletter
\newcommand{\po}{\ext@arrow 0359\rightarrowfill@{}{O}}
\makeatother

\newcommand{\mcf}{Mis\`ere Connect Four\xspace}
\newcommand{\mck}{Mis\`ere Connect $k$\xspace}

\newtheorem{theorem}{Theorem}
\newtheorem{corollary}{Corollary}
\newtheorem{lemma}{Lemma}
\usepackage{amsthm}

\newcommand{\rev}[1]{\textcolor{black}{#1}}

\begin{document}

\author{Robert Steele}
	\email{steele.j.robert@gmail.com}
	\affiliation{NYU Langone Health, New York, NY 10016, USA}
\author{Daniel B. Larremore}
	\email{daniel.larremore@colorado.edu}
	\affiliation{Department of Computer Science, University of Colorado Boulder, Boulder, CO, USA}
	\affiliation{BioFrontiers Institute, University of Colorado Boulder, Boulder, CO USA}
	\affiliation{Santa Fe Institute, Santa Fe, NM, USA}

\title{\mcf is Solved}
\begin{abstract}
	Connect Four is a two-player game where each player attempts to be the first to create a sequence of four of their pieces, arranged horizontally, vertically, or diagonally, by dropping pieces into the columns of a grid of width seven and height six, in alternating turns. \mcf is played by the same rules, but with the opposite objective: {\it do not} connect four. This paper announces that \mcf is solved: perfect play by both sides leads to a second-player win. More generally, this paper also announces that Mis\`ere Connect $k$ played on a $w \times h$ board is also solved, but the outcome depends on the game's parameters $k$, $w$, and $h$, and may be a first-player win, a second-player win, or a draw. These results are constructive, meaning that we provide explicit strategies, thus enabling readers to impress their friends and foes alike with provably optimal play in the mis\`ere form of a table-top game for children.
\end{abstract}
\maketitle

\section{Introduction}

To solve a game is to predict its outcome with certainty when both players play perfectly. For instance, Tic-Tac-Toe, Connect Four~\citep{allis1988knowledge}, and Checkers~\citep{schaeffer2007checkers} are solved, while Chess and Go are not. This means that, although there exist computer programs that can beat any living human in Chess or Go, it is unknown what {\it perfect} play consists of, or how a perfect game ends. In contrast, we know with certainty that perfect play leads Tic-Tac-Toe and Checkers to a tie, while Connect Four ends with a first player win.

The solutions to Tic-Tac-Toe, Connect Four, and Checkers are all rather different from each other. Tic-Tac-Toe has such a small game tree that one can write out each player's optimal move for any opponent move on a single sheet of paper and show that both player 1 (P1) and player 2 (P2) can avoid losing via perfect play. Moreover, this complete solution can still guide one player optimally even against an imperfect opponent.
In contrast, the game tree of Connect Four is far too large to write out, so the first published solution to Connect Four~\citep{allis1988knowledge} instead introduced a set of rules and proved that P1's application of them guarantees a win.\footnote{In fact, in 1988, both James D.\ Allen and Victor Allis solved Connect Four, independently, within 15 days of each other. They had never heard of each other until Allen posted his announcement on 1 October. Allis's announcement followed 15 days later. Allen's proof is not written up as a paper, but he does have a wonderful book on Connect Four~\citep{allen2010connect4}.} However, the existence of these rules---a so-called knowledge-based approach---does not mean that one can use them from memory to beat an {\it imperfect} opponent.\footnote{A note on Allen's \href{https://fabpedigree.com/james/connect4.htm}{webpage} speaks clearly to the difference between perfect algorithmic play and expert human play: ``Eventually I proved by computer that the first player can force a win, but paradoxically when I play against an experienced human, I almost always win moving second, while the opponent sometimes salvages a draw when I move first.''} 
Moreover, for Checkers, perfect play is out of reach of humans entirely, as the solution to Checkers is purely computational, and not rule-based. It required three computer programs to intelligently search from the game's initial position forward, combined with a painstakingly computed endgame database of ten trillion solutions for 10 pieces or fewer which enabled search backward~\citep{schaeffer2007checkers}. A human can implement optimal Tic-Tac-Toe against any opponent, and can understand perfect Connect Four against a perfect opponent, but truly perfect Checkers is played only by machines.

Solutions to games can also be categorized along a different axis, proposed by Allis~\citep{allis1994searching}: an {\it ultra-weak solution} proves whether a player will win, lose, or draw from the game's initial state under perfect play by both sides. A {\it weak solution} provides an algorithm for perfect play by one player against any possible play by the opponent, starting from the game's initial position. And a {\it strong solution} provides an algorithm for perfect play for both players from any legal position. In other words, an ultra-weak solution proves the outcome, a weak solution shows how to achieve it, and a strong solution comprises a weak solution for every starting condition. Through this lens, the 1988 results from Allen and Allis were weak solutions to Connect Four, which were followed by a strong solution due to John Tromp around 1995~\citep{tromp1995website,tromp1995data}. \rev{While weak solutions are intellectually stimulating, strong solutions can be staggering in scale: the number of possible board positions on a typical $7 \times 6$ Connect Four board was computed in 2008 to be $4,531,985,219,092$~\citep{edelkamp2008symbolic}}.

One interesting type of game for theorists and players alike is the {\it mis\`ere game}. A mis\`ere game is a game played by typical rules, but with the winning condition reversed. For instance, in Mis\`ere Tic-Tac-Toe, the loser is the first to get three in a row. In Mis\`ere Checkers, after the addition of a rule making capturing compulsory, the winner is the first to have no pieces remaining. Thus, in \mcf, one's goal is to {\it avoid} connecting four. 

\begin{table*}[t]
\centering
	\begin{tabular}{|c|c|c|c|}
		\multicolumn{4}{c}{Mis\`ere Connect $2$}\\
		\hline
 		\diagbox[]{{\bf height}}{{\bf width}} & $\mathbf{w\geq3}$ {\bf odd}& $\mathbf{2}$, $\mathbf{4}$, {\bf or} $\mathbf{6}$ & $\mathbf{w\geq8}$ {\bf even} \\
		\hline 
		$\mathbf{1}$ & P2 & Draw & P2\\
		$\mathbf{h\geq2}$ & P2 & P2 & P2 \\
		\hline
		\multicolumn{4}{c}{\ }
	\end{tabular}
	\quad
		\begin{tabular}{|c|c|c|}
		\multicolumn{3}{c}{Mis\`ere Connect $k\geq3$}\\
		\hline
 		\diagbox[]{{\bf height}}{{\bf width}} & $\mathbf{w\geq k}$ {\bf odd} & $\mathbf{w\geq k}$ {\bf even} \\
		\hline 
		$\mathbf{1}$ & Draw & Draw \\
		$\mathbf{h\geq2}$ {\bf even} & P2 & P2 \\
		$\mathbf{h\geq3}$ {\bf odd}  & P1 & P2\\
		\hline
	\end{tabular}
	\caption{Outcomes under optimal play for all parameters of \mck. If $w$ or $h$ is infinite, or $w<k$, optimal play leads to a Draw.}
	\label{table-summary}
\end{table*}

Mis\`ere games are practically, mathematically, and computationally interesting. 
Practically, optimal play in mis\`ere games often involves moves and strategies that would be considered poor or counterproductive in normal play, which can be novel, fun, and challenging to experienced players---one cannot simply ``reverse'' one's tactics and strategy. 
Mathematically, analysis of mis\`ere games has challenged combinatorial game theorists because some of the common techniques and arguments used to analyze normal games do not apply to their mis\`ere forms (but see Milley and Renault's review of recent progress~\citep{milley2018restricted}). Indeed, one result has even shown that there is no relationship between normal and mis\`ere outcomes, such that for every pair of (not necessarily distinct) outcomes $\mathcal{O}_1, \mathcal{O}_2 \in \left\{P1, P2\right\}$, there exists a game with normal outcome $\mathcal{O}_1$ and mis\`ere outcome $\mathcal{O}_2$~\citep{mesdal2007simplification}.\footnote{In fact, the result is even more general, such that the set of outcomes includes not only P1 and P2 wins, but also so-called Left and Right wins~\citep{mesdal2007simplification}, in which the ``Left'' or ``Right'' player always wins, respectively, regardless of which goes first. However, this type of outcome does not apply to Connect Four.}
And computationally, solutions to normal games do not imply solutions to their mis\`ere forms, nor vice versa: 
Checkers is weakly solved (draw, \citep{schaeffer2007checkers}), while Mis\`ere Checkers is unsolved.\footnote{So-called {\it mate-in-1} problems in normal-play Checkers can be computed in linear time (in the number of pieces)~\citep{fraenkel1978complexity}, while the mis\`ere problem of {\it lose-in-1} Checkers is NP-complete~\citep{bosboom2019losing}.}
Chess remains unsolved, while Mis\`ere Chess is weakly solved (P1 win, \citep{watkins2017losing}). 

This paper announces that \mcf is weakly solved: perfect play by both sides results in a P2 win. More generally, this paper also announces that \mck played on a board of width $w$ and height $h$ is also weakly solved, but the outcome depends on the games parameters $k$, $w$, and $h$, and may be a P1 win, a P2 win, or a draw. 
Finally, we solve ``infinite'' \mck by addressing boards of infinite height or infinite width, both of which are draws. A summary of these results is given in Table~\ref{table-summary}.

\section{Rules and Notation}

\mck is played on a board of height $h$ and width $w$ which forms a grid. Player 1 (P1) plays as \verb |X| and Player 2 (P2) plays as \verb |O|. In alternating turns, each player selects a column in which to place a piece, and gravity pulls that piece to the lowest open position. Consequently, once a column has been played $h$ times, it is full and neither player can play in it. The game ends when one player has $k$ pieces in a contiguous horizontal, vertical, or diagonal line. Because the game is played in the mis\`ere condition, the first player with $k$ pieces in a row loses. 

Using common vocabulary from the literature, this makes \mck is a two-player partisan\footnote{Connect Four is {\it partisan} because it is not {\it impartial}. In impartial games games, allowable moves depend only on who is playing. Connect Four, along with Go, Chess, and Checkers, are not impartial because the two players play different colored pieces.} game with perfect information\footnote{Perfect information means that nothing is hidden from anyone.}, played under the mis\`ere condition: the player who connects four loses. 

\section{Solving the game}

Our solution is divided into six parts.
First, we introduce a strategy and prove that the second player (P2) can always use it to win \mcf on a standard $6\times7$ board. 
Second, using the same strategy, we extend the proof to show that P2 wins \mck on all boards of even height $h$ and $k\leq w$.
Third, we extend the strategy to cover boards of odd height $h\geq3$, proving that P2 wins when the width $w\geq4$ is even, but P1 wins when $w\geq5$ is odd. 
Fourth, we analyze games of height $h=1$ for $k\geq 3$. Fifth, we analyze games for $k=2$. 
Finally, we treat two special cases: narrow boards $w<k$ and infinite boards.

Perhaps surprisingly, the proofs for $k=2,3$ and $h=1$ turn out to be the most technical and important results in the paper. In fact, the proof for $k=3$ and $h=1$ underpins some of the solutions in the sections that precede it. Due to that importance, we state its result informally here: {\it On a single-row board (height $h=1$) of arbitrary width, neither player can be forced to connect $k\geq3$, {\it even if} we allow their opponent to connect as many as they wish.} This key result means that all $h=1$ and $k\geq3$ games are a draw.

\subsection{\mcf on a standard $6\times 7$ board}

\begin{figure}[t]
	\centering
	\includegraphics[width=1.0\linewidth]{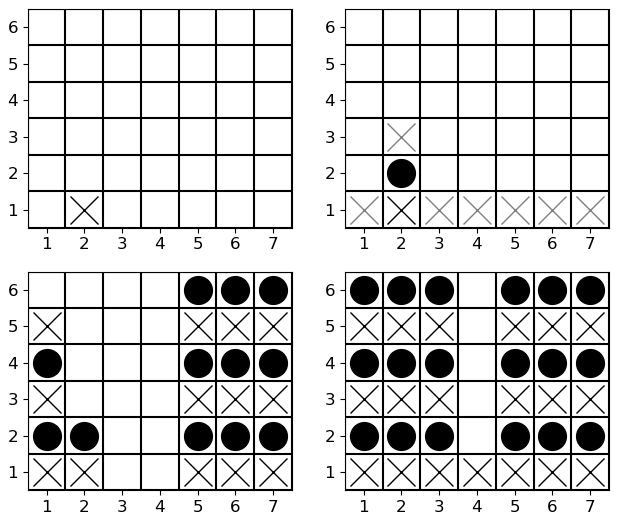}
	\cprotect\caption{{\bf Illustration of the \textit{take-even} strategy}. (Top-left) In the first move of the game, P1 (\verb|X|) must play into row 1. (Top-right) The take-even strategy dictates that P2 (\verb|O|) plays in an even row. The existence of this move is guaranteed by P1's play in an odd row. Grey \verb|X|'s indicate the available next moves for P1, all of which are in odd rows. (Bottom-left) An example game state after 25 turns, in which P2 is using the take-even strategy, with P2 to play. (Bottom-right) Example final game state with a P2 win on turn 37. This solution is proved in Theorem~\ref{thm-take-even}.}
	\label{fig-take-even}
\end{figure}

The second player (P2) can always win \mcf played on a standard board of height $h=6$ and width $w=7$. P2 achieves this victory using a simple strategy: playing only in the even rows. As a consequence P1 is eventually forced to connect in the first row, resulting in a P2 win. We prove that P2's {\it take-even} strategy results in an unavoidable P2 win in Theorem~\ref{thm-take-even}, and illustrate it in Figure~\ref{fig-take-even}.  

\begin{theorem}\label{thm-take-even}
	The take-even strategy, played by P2, on a board with height $h=6$ and width $w=7$, eventually forces P1 to connect four resulting in a guaranteed P2 win.
\end{theorem}
\begin{proof}
When the game begins, the board is empty so P1 has no choice but to play in the first row---an odd-numbered row. By playing in an odd row, P1 has created an opportunity for P2 to play in an even row. In fact, the space on top of P1 is the only even-row move that P2 has available. According to the take-even strategy, P2 plays on top of P1. This fills the even-row space, confronting P1 with only odd-row options. 
	
This opening illustrates three important facts. First, if P1 has only odd-row plays available, then P1 must play in an odd row. Second, if P1 plays in an odd row, this creates an even-row play for P2 on top of P1's most recent piece. P2's even-row option is guaranteed, in fact, because the board is even height, so P1's odd-row play must have at least one empty space above it. Third, when P2 takes the even-row play that P1 has created, P2 removes the even-row option from the board, confronting P1 with only odd-row options. Note that if P2's move completely fills a column, the number of playable columns decreases but P1 has only odd-row options. 

Together, these facts mean that P2 always has the option to play in an even row and that P1 is forced to play in odd rows. In other words, P2 can exercise the take-even strategy, and P1 must play in odd rows. P2 has control of the game.

Now, note that because P2 plays only in even rows and P1 plays only in odd rows, neither player can connect four vertically or diagonally. Instead, P1 will eventually be forced to connect four in the first row, once the other columns are filled (by P2) and thus unavailable. P2 therefore wins on the 37th move, or sooner.
\end{proof}

\subsection{Mis\`ere Connect $k\leq w$ on a board with even height}

The take-even strategy works on the standard $6\times7$ board because of three key ingredients. First, the board's even height guarantees that any odd-row play by P1 will have an empty even-row space above it, into which P2 will play. Second, the board is finite so eventually the columns fill forcing P1 to play in the first row. Third, the board is wide enough that if P1 fills the first row, P1 loses, i.e. $k \leq w$. The fact that these key ingredients are not restricted to $6 \times 7$ boards alone leads to the following, more general Corollary to Theorem~\ref{thm-take-even}.

\begin{corollary}\label{cor-even-boards}
	The take-even strategy played by P2, on a board with even height $h$ and finite width $w$, with $w \geq k$, eventually forces P1 to connect $k$ resulting in a guaranteed P2 win.
\end{corollary}
\begin{proof}
	The proof of Theorem~\ref{thm-take-even} applies directly without modification. P1 can delay loss as long as possible by playing only in sets of $k-1$ adjacent columns. Under this delay tactic, with $w$ total columns and at most $\lfloor \tfrac{w}{k}\rfloor$ individual columns may be empty. Each column allows $h$ moves. Thus, P1 loses on or before move number $1+\left(w-\lfloor \tfrac{w}{k}\rfloor\right)h$.
\end{proof}

\subsection{Mis\`ere Connect $k\leq w$ on a board with odd height $h\geq 3$}

An even board height $h$ is required by the assumptions of Theorem~\ref{thm-take-even} and  Corollary~\ref{cor-even-boards}, but what if the board's height is odd? At first glance, an odd height appears to undermine a key pillar of the take-even strategy, namely the fact that any odd-row move by P1 creates an even-row option for P2. Now, with an odd number of rows, P1 could play in the top row, filling a column but creating no even-row option for P2's strategy. Without an even-row option, P2's take-even strategy, by definition, falls apart.

While at first this appears to spell trouble for P2 (and our proof strategies), we now show that, in fact, the machinery of Theorem~\ref{thm-take-even} can still be used for odd-height boards, but with a twist: the player whose win is guaranteed now depends on the {\it width} of the board, with P2 winning even-width games but P1 winning odd-width games.

\begin{figure}[t]
		\centering
		\includegraphics[width=0.8\linewidth]{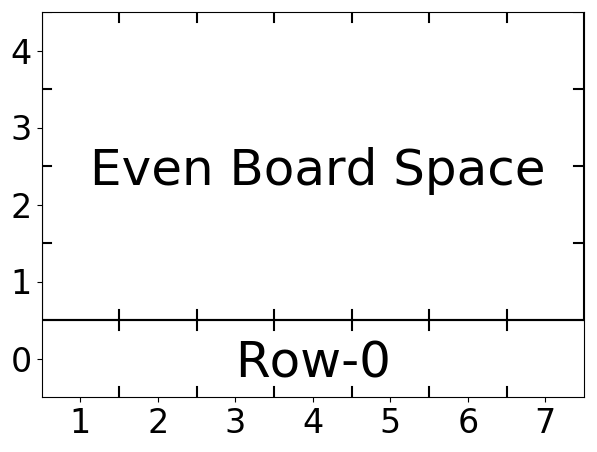}
		\caption{The delayed take-even strategy works by conceptually partitioning the board. The bottom row is counted as "row zero" and the rows above, of which there must be an even number, are called the ``even board space.''}
		\label{fig:even_board_space}
\end{figure}

To build intuition toward our next Theorem, consider a board with odd height where we have numbered the bottom row $0$, and the rows above as $1, 2, \dots, h-1$. Because $h$ is odd by assumption, $h-1$ is even, and so the board can be divided into two sections: row $0$, and a section on top with an even number of rows. We call the area on top of row $0$ with an even $h-1$ number of rows the {\it even board space} (Fig.~\ref{fig:even_board_space}).

With this partition of the board in hand, we introduce the {\it delayed take-even} strategy for boards of odd height $h$---an unbeatable strategy for P2 when width $w$ is even, and for P1 when width $w$ is odd. Here is how it works, in two simple rules: (1) if the opponent plays in row $0$, follow with a non-losing move in row $0$; (2) if the opponent plays in any other row, play on top of them, thus using a take-even strategy in the even board space. 

The logic of this strategy is simple: if row $0$ were to be filled, then the game could be played as if there were an even number of rows, thus allowing the take-even strategy to guarantee a win. Indeed, the take-even strategy for even boards tells us that the first player to play will lose. In odd height boards the finding is similar: the first player to play {\it in the even board space} loses. When there are an even number of moves available in row $0$, P2 can force P1 to be the first to enter the even board space (and thus lose); when there are an odd number of moves available in row $0$, P1 can force P2 to be the first to enter (and thus lose). The following theorem proves this optimal strategy for $k \geq 3$ and $w \geq k$.

\begin{theorem}\label{thm-delayed}
	The delayed take-even strategy, played by P2 on a board with odd height $h\geq3$ and even width $w \geq k \geq 3$, eventually forces P1 to connect $k$ resulting in a guaranteed P2 win.
\end{theorem}
\begin{proof}
	P1 has no choice but to begin the game by playing in row $0$. Because P1 has played in row $0$, the delayed take-even strategy dictates that P2 should also play in row $0$. No one has yet lost (or won) because each has played just once and $k\geq3$ by assumption. Thus, P1 now has a choice between two options: (i) play into an odd-numbered row (that is, row $1$, in the even board space, on top of one of the existing moves), or (ii) play again into row $0$. 
	
	Note that no matter which option P1 chooses, P2 will return to P1 a board with only option (i), option (ii), or both, meaning that P1 will not have the option to play into an even-numbered row in the even board space. Why must this be so? First, if P1 plays into an odd row, P2 will play on top, either giving P1 another odd-row option on top, or completing a column and removing it from P1's options entirely. The availability of P2's move is guaranteed by the fact that there are an even number of rows in the even board space above row $0$.
	On the other hand, if P1 plays into row $0$, P2 will also play a non-losing move into row $0$, either giving P1 another row-$0$ option or completing row $0$.  The availability of P2's move is guaranteed by the fact that there are an even number of {\it columns} in row $0$, so that any time P1 plays in row $0$, there is at least one column into which P2 may play. Critically, the guarantee that such a move by P2 is non-losing is established in Theorem~\ref{theorem-pancake} which shows that both players can bring play in row $0$ to a draw with perfect play. 
	
	Taken together, when P2 uses the delayed take-even strategy, P1's options are restricted entirely to playing in row $0$ or playing in an odd row. And, regardless of P1's choice, P2 will have an available move that returns P1 to these identical choices. Eventually, if P2 has not already won in the even board space by Theorem~\ref{thm-take-even} and Corollary~\ref{cor-even-boards}, row $0$ will be filled, and the game will be {\it forced} into the even board space, leading to a P2 win.
\end{proof}

The proof of Theorem~\ref{thm-delayed} shows how the delayed take-even strategy can be utilized by P2 to eventually force P1 into the even board space leading to a P2 win, provided that the height of the board is odd and the width of the board is {\it even}. Using similar arguments, we now show that when the height of the board is odd and the width of the board is {\it odd}, the delayed take-even strategy reverses the outcome, allowing P1 to force P2 into the even board space leading to a P1 win. 

\begin{corollary}\label{cor-delayed}
The delayed take-even strategy, played by P1 on a board with odd height $h\geq3$ and odd width $w \geq k \geq 3$ eventually forces P2 to connect $k$ resulting in a guaranteed P1 win.
\end{corollary}
\begin{proof}
	P1 has no choice but to begin the game by playing row $0$. P2 now has a choice between two options: (i) play into an odd row (in the even board space, on top of one of the existing moves), or (ii) play again into row $0$. Identical arguments to Theorem~\ref{thm-delayed} now apply: P1's use of the delayed take-even strategy will always return P2 to these two choices. However, in contrast with Theorem~\ref{thm-delayed}, this Corollary assumes the width to be odd, allowing P1 to play the final move in row 0. Thus, P1 can force P2 into the even board space, eventually resulting in a P1 win. Once more, this result rests on the fact that P1 can avoid a loss during gameplay in row 0, a result which we establish in Theorem~\ref{theorem-pancake}.
\end{proof}

Theorem~\ref{thm-delayed} and Corollary~\ref{cor-delayed} prove all cases of $k\geq 3$ for both odd and even width, and both odd and even height, provided that there is sufficient width to connect $k$ in a row, i.e., $w\geq k$. However, as noted, both proofs rely on an as-yet unproven Theorem, one which guarantees that P1 always has a non-losing move available in row $0$ for odd-width boards, and that P2 always has a non-losing move available in row $0$ for even width boards. This result follows from proving that perfect play by both players results in a draw for boards with height $h=1$ when $k\geq3$, the precise result of Theorem~\ref{theorem-pancake}.

\rev{\subsection{\mck with $k\geq3$ and height $h=1$}}

When a board has height $h=1$, \mck simply asks whether it is possible for either player to force the other to occupy $k$ adjacent positions along a strip of some width. The answer, as we will prove, is that it not possible: \mck with $h=1$ and $k\geq 3$ is a draw. In fact, a stronger statement can be made: neither player can force the other to connect even under a surprising relaxation of the rules! For $k=3$, either player can avoid connecting $3$, {\it even if} we allow the other player to connect as many pieces as desired without losing. In other words---and perhaps to our surprise---P1 can draw the $h=1$ game even if P2 cannot lose, and P2 can draw the $h=1$ game even if P1 cannot lose.

The definitions and strategies below, we let P1 play as \verb|X| and P2 as \verb|O|. The core of the proof relies on enumerating all the empty spaces on the board at the start of the game, in pairs, from left to right.\footnote{This enumeration and pairing argument should be attributed to one of the paper's anonymous reviewers, to whom we are grateful. Prior to their suggestion, our proof was substantially longer, more complicated, and less elegant. That earlier proof of our claims can be found in Appendix~\ref{alternative}, and relies on analyzing the board as a set of automata which either player can always manage toward a draw. See} For instance, for a game of $w=8$ or $w=11$, we would label the spaces
$$\verb|   board:  --------|\vspace{-0.3in}$$
$$\verb|  labels:  11223344|$$
\vspace{-0.3in}
$$\verb|   board:  -----------|\vspace{-0.3in}$$
$$\verb|  labels:  11223344556|$$
If the game ends in a balanced state such that every pair has one \verb|X| and one \verb|O|, the game is a draw for any $k\geq 3$, as there can be at most two-in-a-row. Note that this is true regardless of who has played in the unpaired rightmost space when $w$ is odd. We refer to this as a {\it pair-balanced} board.

Now consider a simple two-rule strategy, which can be employed by either player, and is designed to force any game to a pair-balanced board, and thus, a draw. The two rules are:
\begin{enumerate}
	\item Take-other: If your opponent has taken one half of a pair, take the other half of the pair. 
	\item Rightmost: If there are no available take-other moves, play the rightmost space available.
\end{enumerate}
In the theorem that follows, we prove that this two-rule strategy leads to a draw because it guarantees that the board ends in a pair-balanced state, either (i) by showing that every pair is balanced, or (ii) by showing that no pair is \verb|XX| or \verb|OO|.

\begin{theorem}\label{theorem-pancake}
	Suppose that $h=1$ and $k\geq 3$. Either P1 playing \verb|X| or P2 playing \verb|O| can force a draw by playing the two rules of take-other and rightmost.
\end{theorem}
\begin{proof}
	Case 1: Suppose that $w$ is even and \verb|O| is employing the two-rule strategy. By definition, \verb|X| plays somewhere on the board, half filling a pair. Because \verb|O| uses take-other, a balanced pair is created. With only empty pairs to play in, \verb|X|'s subsequent play must half fill another pair, which \verb|O| will then complete. Because $w$ is even, the game consists of $w/2$ repetitions of this process, and the final board is guaranteed to be pair-balanced.
	
	Case 2: Suppose that $w$ is even and \verb|X| is employing the two-rule strategy. Their first move fills half of the rightmost pair, creating a \verb|-X|. When \verb|O| plays, they will either fill that \verb|-X| or half fill some other pair. If \verb|O| fills the \verb|-X|, that pair is balanced and the game is equivalent to a new $w-2$ game. If \verb|O| half fills some other pair, then \verb|X| plays take-other to balance that pair, and once more, \verb|O| has two options: balance the \verb|-X| or half-fill another pair. Thus, \verb|X| always has an in-strategy option and is never forced to play the other half of the \verb|-X|. Because $w$ is even, \verb|O| plays last, and the final board is guaranteed to be pair-balanced.
	
	Case 3: Suppose that $w$ is odd and \verb|O| is employing the two-rule strategy. Because $w$ is odd, the first move from \verb|X| either takes the rightmost space, which is not part of any pair, or \verb|X| half fills a pair elsewhere on the board. Suppose that \verb|X| takes the rightmost space. In this case, it is as if \verb|O| is playing first on a $w-1$ board, and so the arguments of Case 2 (with the players reversed) guarantee a pair-balanced final board. If instead, \verb|X| takes some other space, \verb|O| will take-other to balance the pair, and we restart Case 3 with $w-2$. 
	
	Case 4: Suppose that $w$ is odd and \verb|X| is employing the two-rule strategy. By the rightmost rule, \verb|X| takes the singleton space, leaving \verb|O| to play somewhere in the even $w-1$ board. The arguments of Case 1 (with the players reversed) guarantee a pair-balanced final board.
	
	These four cases cover both even and odd width $w$ for both players, meaning that either player can force $h=1$ and $k\geq 3$ to a draw.
\end{proof}

With the result of Theorem~\ref{theorem-pancake}, the results of Theorem~\ref{thm-delayed} and Corollary~\ref{cor-delayed} are complete: because \mck with $h=1$ and $k\geq 3$ is a draw, neither player can win in the bottom row (row 0) of an odd-height board, which means that the delayed take-even strategy cannot be escaped: for odd heights $h\geq 3$, P2 wins when width $w$ is even, and P1 wins when width $w$ is odd.

\subsection{\mck with $k=2$}

When $k=2$, the game takes on a fundamentally different character. Connecting two requires only that a player place two pieces adjacent or diagonal to each other, leading to the question: can one player force another to play adjacent to one of their previous moves? The answer, as we shall see, is generally yes, and in particular, the advantage belongs to P2.

\rev{We begin by analyzing the special case of $k=2$ with height $h=1$ which will be subsequently used to solve $k=2$ and arbitrary $h$.} Playing \mck with $k=2$ with height $h=1$ amounts to a game in which each player must simply avoid placing two pieces adjacent to each other. This simple goal leads to a few preliminary and useful observations. 

First, the game is a draw if and only if it ends with alternating \verb|X| and \verb|O| pieces. This means that there are only two types of final boards for any game that ends in a draw: one in which \verb|X| has placed a piece in the leftmost position, and one in which \verb|O| has. Consequently, once \verb|X| has played the first piece, the final game state that would constitute a draw is uniquely specified. We call this the board's {\it template}. A move that complies with the template is an ``in-template'' move, and a move that does not---and thus asserts a contradicting template of its own---is an ``anti-template'' move. 

If either player chooses to play an anti-template move after the opening move has asserted a template, then the game must have a winner and a loser. Thus, the key question to be addressed is whether playing an anti-template move is strategically useful in allowing one player to force the other to connect 2. 

We now introduce a simple bookkeeping scheme, which tracks the number of non-losing in-template moves available to \verb|X| and \verb|O|, respectively, as an ordered pair. For instance, the following board has 1 available in-template move for \verb|X| and 2 for \verb|O|, as annotated.
$$\verb|   board:  -XO--|\vspace{-0.3in}$$
$$\verb|           ↑  ↑↑|\vspace{-0.3in}$$
$$\verb|template:  O  XO|$$
We would therefore write the number of available in-template moves at $(1,2)$.

This bookkeeping scheme extends to boards with two contradictory templates. For instance, we can count $(1,1)$ available template moves in the board \verb|X---O| if we assert \verb|X|'s template and  annotate from left to right:
$$\verb|   board:  X---O|\vspace{-0.3in}$$
$$\verb|           ↑↑ |\vspace{-0.3in}$$
$$\verb|template:   OX  |$$
Note that we count an identical number if we assert \verb|O|'s template annotate from right to left:
$$\verb|   board:  X---O|\vspace{-0.3in}$$
$$\verb|            ↑↑|\vspace{-0.3in}$$
$$\verb|template:    OX |$$
In fact, we can assert both templates and annotate from the outside in---we again get an identical number:
$$\verb|   board:  X---O|\vspace{-0.3in}$$
$$\verb|           ↑ ↑|\vspace{-0.3in}$$
$$\verb|template:  O X|$$
This example illustrates the point that we can choose any template as our reference to tabulate the number of remaining moves for each player, and the result will be the same. It is the contradiction between two templates that reduces the number of available moves; not the order of tabulation. In some sense, our bookkeeping scheme simply imagines that both players play alternating non-losing moves arbitrarily until there are none left.

Our bookkeeping scheme is useful because it tells us what will happen if both players play alternating non-losing moves. For instance, if there exist both in- and anti-template moves on the board (guaranteeing a winner and a loser), and the move tally is $(1,0)$ with \verb|X| to play, then \verb|X| will win: \verb|X| has an available move and will play it, while \verb|O| will have no moves available. By extension, if the tally is $(a+1,a)$ with \verb|X| to play, then \verb|X| will win, unless of course \verb|O| plays outside the existing templates to intentionally shift the tally. We will now explore such tally-shifting plays by both players.  

From the perspective of the tally, there are only four canonical plays. Without loss of generality, we'll define them from the perspective of \verb|X|, though they equally apply to \verb|O|. 

The first canonical play is the {\bf in-template play}, in which \verb|X| plays into an existing template. For instance, 
$$\verb|X----O| \px \verb|X-X--O|$$
is an in-template play. The tally keeps track of available in-template plays, so tabulating it before and after \verb|X|'s play shows $(2,2) \px (1,2)$, as expected. The in-template play affects only the person who played it, decreasing their tally by one. We therefore give the in-template play a score of $(-1,0)$ for its effect on the tally.

Note that an in-template does not require that all plays on the board conform to one template. For instance, returning to a previous example, in which two contradictory templates have already been asserted, note that
$$\verb|X---O| \px \verb|X--XO|$$
is also an in-template play, such that $(1,1) \px (0,1)$. Thus, this play requires only that \verb|X| play in {\it one of} the templates. Any play that does not create a new contradiction with an existing template is an in-template play.

The second canonical play is the {\bf double contradiction play}, in which \verb|X| plays in contradiction to an existing template that is bounded on both the left and the right by existing moves, thus creating a first contradiction to the left of the play and a second contradiction to the right. For instance,
$$\verb|X------O| \px \verb|X--X---O|\vspace{-0.3in}$$
$$\verb| OXOXOX | \hspace{0.26in} \verb| O- -OX |$$
shows a play (top line), its template before the play (bottom left) and one possible templating afterward (bottom right). Note that this play created not one but two contradictions. Consequently, the play itself debits once from the tally, and the two contradictions debit once each. Tabulating the tally before and after \verb|X|'s play shows $(3,3) \px (1,2)$. We therefore give the double-contradiction play a score of $(-2,-1)$ for its effect on the tally. Note that a double contradiction requires a minimum of three empty spaces between the template walls (a consequence of the fact that there must be at least three template moves that can be depleted from the tally).

The third and fourth canonical plays are when \verb|X| plays in contradiction to an existing template that is bounded on only one side by an existing move, and meets the edge of the board on the other. Both plays are played in the same manner, but their effects differ dramatically, depending on which player is prescribed the wall-adjacent place in the existing template. Consider the following two examples. 

In the first example, \verb|X| {\it is not} prescribed the wall in the existing template, 
$$\verb|XO------| \px \verb|XO---X--|\vspace{-0.3in}$$
$$\verb|  XOXOXO| \hspace{0.26in} \verb|  -XO OX|$$
so when \verb|X| plays a contradiction, it amounts to $(3,3) \px (2,2)$, thus scoring $(-1,-1)$. We call this play an {\bf offensive play} because it decreases \verb|X|'s tally by one like an in-template play, but also decreases \verb|O|'s tally. An offensive play is strictly better than an in-template play.

In this first example, one might also note that \verb|X| left the wall position available after flipping its template. Consequently, \verb|O| has an opportunity to counter with an offensive play of their own, reversing \verb|X|'s short-lived advantage. Therefore, we note that an offensive play directly against the wall is advised, as it prevents one's opponent from countering with a followup offensive play. We call the offensive play against the wall an {\bf exclusive offensive play} for this reason. 

In a second example, \verb|X| {\it is} prescribed the wall in the existing template, 
$$\verb|XO-----| \px \verb|XO---X-|\vspace{-0.3in}$$
$$\verb|  XOXOX| \hspace{0.26in} \verb|  -XO O|$$
so when \verb|X| plays a contradiction, it amounts to $(3,2) \to (1,2)$, thus scoring $(-2,0)$. We call this play a {\bf self-immolation} because it depletes two from \verb|X|'s tally, but none from \verb|O|'s. A self-immolation is strictly worse than an in-template play. 

\begin{table}[h]
	\renewcommand{\arraystretch}{1}
	\begin{tabular}{| l | c | c |}
		\hline
		play name & $\Delta$ tally (\verb|X| played) &  $\Delta$ tally (\verb|O| played)\\
		\hline
		in-template & $(-1,0)$ & $(0,-1)$ \\
		double contradiction & $(-2,-1)$ & $(-1,-2)$ \\
		offensive & $(-1,-1)$ & $(-1,-1)$ \\
		self-immolation & $(-2,0)$ & $(0,-2)$\\
		\hline
	\end{tabular}
	\cprotect\caption{{\bf The four canonical plays for \mck with $k=2$ and $h=1$.} The four named plays, being described in the text, are tabulated here with their effects on the tally of remaining in-template moves. One column shows the impact of such a play when \verb|X| makes that play, and the other column shows the impact when \verb|O| makes that play.}
	\label{table-contradictions}
\end{table}

These plays now allow us to analyze gameplay more generically. For one player to beat another, someone must create a contradiction, and only double contradiction, offensive, or self-immolation moves do so. Relative to in-template moves, double contradiction plays affect both players equally, offensive plays improve the standing of the one who plays, and self-immolations are bad for their player.

\begin{lemma}\label{lemma-k2-odd}
	In a game of \mck with $k=2$, $h=1$, and odd width $w\geq3$, Player 2 wins. 
\end{lemma}
\begin{proof}
	If $w\geq3$ is odd, then $w=2m+1$ for a natural number $m\geq1$. The first move by P1 (i.e., \verb|X|) will either template the board to have $(m-1, m+1)$ template spaces remaining for each player or $(m,m)$ template spaces remaining for each player. The latter outcome happens when \verb|X| plays into an odd-numbered space (indexing from $1$ to $w$, left to right), and the former outcome happens when \verb|X| plays into an even-numbered space.

First, suppose that \verb|X|'s opening move is in an even-numbered space, leaving the tally at $(m-1,m+1)$ with \verb|O| to play. If both players play in-template moves for the rest of the game, \verb|O| wins by tally. Thus, the \rev{burden} will be on \verb|X| to rescue themself from the opening play. Without yet considering \verb|O|'s move, note that \verb|X|'s only path to victory is if \verb|X| is able to play two offensive moves to gain the advantage in the tally. Unfortunately, \verb|O| plays second, and can therefore take one of the in-template wall positions ($1$ or $w$), removing any possibility of an offensive move by \verb|X| on that side. With at most one offensive move available to \verb|X|, \verb|O| can win by tally by playing only in-template moves. 

Instead, suppose that \verb|X|'s opening move is in an odd-numbered space, leaving the tally at $(m,m)$ with \verb|O| to play. Note that \verb|X| must have left either one or both walls open and templated as \verb|X|. Thus, \verb|O| can play an exclusive offensive move by playing an available wall. The tally is now guaranteed to be $(m-1, m-1)$ with \verb|X| to play, and a contradiction on the board. Because there are no offensive moves available to \verb|X|, \verb|X| has no opportunity to reverse their fate, and thus \verb|O| can win by tally by playing only in-template moves. 

In short, if \verb|X| opens an odd-width game in an even numbered space, \verb|O| should play an in-template move against a wall. If \verb|X| opens an odd-width game in an odd numbered space, \verb|O| should play an offensive move against an open wall. In every possible scenario for odd width $w\geq3$, $h=1$, and $k=2$, \verb|O| wins. 
\end{proof}

Even-width games are more interesting. Starting from $w=4$, we note that if \verb|X| opens by playing against a wall, then \verb|O| should take the opposite wall, guaranteeing a draw. If \verb|X| opens by playing an interior space, \verb|O| plays an offensive move and \verb|X| loses. Thus, \verb|X| should open against a wall, and $w=4$ is a draw. 

Similarly, when $w=6$, \verb|X| should open against a wall to draw the game. Why? If \verb|X| were to open in space 2 or 3, then \verb|O| could play an offensive move and win by tally. This holds true even when \verb|X| is able to play an offensive move of their own, a possibility when \verb|X| opens in space 3. Importantly, the path to a draw with $w=6$ is narrow: after \verb|X| opens against a wall, \verb|O| must take the other wall to prevent \verb|X| from winning. And, after \verb|O| takes that wall, \verb|X| must play in-template closest to her first piece, avoiding three empty spaces in a row into which \verb|O| could play a double contradiction to win. Exhaustively, $w=6$ is a draw.  

Unfortunately for \verb|X|, when $w=8$, \verb|X| loses the ability to draw the game, which we now prove.
\begin{lemma}\label{lemma-k2-even}
	In a game of \mck with $k=2$, $h=1$, and even width $w\geq8$, Player 2 wins. 
\end{lemma}
\begin{proof}
First, assume that P1 (i.e., \verb|X|) opens with a wall play. Then \verb|O| is free to take the opposite wall, an in-template play. The entire board is now templated with zero contradictions, and the tally is $(m-1,m-1)$ with \verb|X| to play. Because the walls are both taken, there are only in-template and double contradiction moves remaining. If \verb|X| plays a double contradiction, the tally will drop to $(m-3,m-2)$, and \verb|X| will lose by tally. Thus, \verb|X| must play only in-template moves. After \verb|X| plays an in-template move, the tally will be $(m-2,m-1)$ with \verb|O| to play. If \verb|O| plays a double contradiction, the tally will drop to $(m-3,m-3)$ with \verb|X| to play, and \verb|X| will lose by tally. The existence of the double contradiction for \verb|O| to play on their second turn is guaranteed: \verb|O| needs three or more empty spaces in which to play this move, and \verb|X| can break the stretch of six spaces into a length of two and a length of three.\footnote{Incidentally, this is why we needed to treat $w=4$ and $w=6$ as special cases, and why their outcomes differ from the $w\geq8$ case! With fewer than $8$ total spaces, the availability of the double contradiction is no longer guaranteed.} Consequently, if \verb|X| opens with a wall play, \verb|O| can win by taking the opposite wall on their first turn and a double-contradiction on their second, with in-template moves thereafter to win by tally. 
	
Next, assume that \verb|X| opens with something other than a wall play. This means that both walls are open, one of which must be an offensive move for \verb|O|. Thus, \verb|O| plays the offensive move, and the tally is $(m-2,m-1)$ with a contradiction on the board. From here, it is hopeless for \verb|X|. Even if \verb|X| has left enough room to play an offensive move against the opposite wall, this will bring the tally to $(m-3,m-2)$ with \verb|O| to play. With no tools to overcome the tally deficit, template moves by \verb|O| will lead to a \verb|O| win.

In short, if \verb|X| opens an even-width $w\geq8$ game against a wall, then \verb|O| should take the opposite wall and then play a double-contradiction. If \verb|X| opens somewhere other than a wall, then \verb|O| should play an exclusive offensive move at the appropriate wall. In every possible scenario for even width $w\geq8$, $h=1$, and $k=2$, \verb|O| wins.
\end{proof}

Our analysis of $k=2$ and $h=1$ games is complete, and points to a general advantage for \verb|O|, who can win for any odd width $w\geq3$ and any even width $w\geq8$. The special cases of $w=4$ and $w=6$ are a draw. We conclude with a final theorem, generalizing these results to arbitrary $h$. 

\begin{theorem}
	In a game of \mck with $k=2$, optimal play leads to a draw for height $h=1$ and widths $w=1$, $2$, $4$, and $6$, and a P2 win for all other finite widths. When $h>1$ and $w=1$, the game is a draw. When $h>1$ and $w>1$, optimal play leads to a P2 win.
\end{theorem}
\begin{proof}
For height $h=1$, widths $w=1$ and $2$ are obviously draws because neither player has a chance to play more than once. All other odd widths are wins for P2 by Lemma~\ref{lemma-k2-odd}. Even widths $w=4$ and $6$ are draws by directly exploring all possible game trees (see text), while all other even widths $w\geq8$ are wins for P2 by Lemma~\ref{lemma-k2-even}. The Lemmata and discussion above provide strategy prescriptions for P2, and thus, this proof is constructive.
	
For height $h>1$ and width $w=1$, the game follows a deterministic trajectory where the players fill the single column in alternative moves until the board is full with no winner. We include this game only for the sake of completion.

For even heights $h\geq2$ and finite width $w>1$, we simply observe that \verb|O| can use the take-even strategy, guaranteeing a P2 win, despite the other peculiarities of the $k=2$ cases.

For odd heights $h\geq3$ and finite width $w>1$, P2 wins. This is because, on each turn, each player must decide whether to play in the bottom row, or to play above in the even board space (Fig.~\ref{fig:even_board_space}). P2 can always win by using a delayed take-even strategy---that is, by following P1's bottom-row moves with bottom-row moves, and otherwise playing on top of P1. This follows from two facts: First, the delayed take-even strategy forces P1 into alternating bottom row play. Second, Lemmata~\ref{lemma-k2-odd} and \ref{lemma-k2-even} provide a strategy by which P2 can either win or draw in alternating bottom-row play. 

One lingering question might be whether, somehow, the existence of additional rows would somehow disrupt the tally arguments used to prove $h=1$ cases. However, we alleviate this concern by noting that (i) for even heights $h$, the take-even strategy drives P2 to victory, not a tally-based strategy, and (ii) for odd heights $h$, the tally-based strategy is never disrupted by the delayed take-even strategy because P2 would never play in row 1, and thus P2 would never decrease their tally of in-template plays in row 0 via the threat of a diagonal connect 2.  
\end{proof}

\subsection{Two special cases: infinite games and narrow games}

Trivially, boards of infinite width must be draws: either player can play arbitrarily far from any previous play. Similarly, boards of infinite height must also be draws: either player can indefinitely play on top of the other's most recent move. Incidentally, infinite Connect Four is also a draw~\citep{yamaguchi2012infinite}.

Finally, if $k$ exceeds the width of the board $w$, the game is a draw: either play can play on top of the other's most recent move, avoiding a vertical connection. Because the board is also too narrow to connect horizontally or diagonally, the game is a draw.

\section{Discussion and Conclusion}

In this paper, we described \mcf, the mis\`ere form of Connect Four in which players must force their opponents to connect four in order to win, and generalized it to \mck played on a board of width $w\geq 1$ and height $h\geq 1$. We proved that under optimal play the game may be a draw, Player 1 win, or Player 2 win, summarized in Table~\ref{table-summary}. In general, the outcome of the game has little to do with $k$, and more to do with the dimensions of the board, particularly whether height and width are even, odd, one, or infinite. 

One observation is that the game gently defies expectations, in that one might expect that the first player to play is at a disadvantage, and thus games ought to end in either a draw or a Player 2 win. However, this is not always the case: when $3 \leq k \leq w$ with $w$ and $h$ both odd, Player 1 wins. 

A second observation is that the solution to \mck is constructive, realizable by a human, and in some sense {\it thick}: over all possible plays by an opponent, the proven winner (or drawer) can apply the knowledge-based strategies of this paper to pull the game toward the optimal outcome. In contrast, when one's opponent plays suboptimally in other games, even a (human) player who understands how a perfectly played game unfolds may struggle to respond optimally. In other words, we have not computed a thread between the opening of an optimally played game and its conclusion, but have instead provided a set of blanket strategies that always pull the game toward its inevitable conclusion.

Finally, we leave one open problem for readers' consideration. Consider \mck played on an annulus, i.e., on a board of height $h$ and circumference $w$ (or width $w$ with periodic walls).  If the height is even, the take-even strategy in this paper can be used without modification, but who wins when the height is odd, including when the height \rev{is} one?

\onecolumngrid
\clearpage
\appendix
\onecolumngrid
\section{Alternative Proof: \mck with $k\geq3$ and height $h=1$}\label{alternative}
\twocolumngrid

When a board has height $h=1$, \mck simply asks whether it is possible for either player to force the other to occupy $k$ adjacent positions along a strip of some width. In fact, it is not possible, even under a surprising relaxation of the rules! For $k=3$, either player can avoid connecting $3$, {\it even if} we allow the other player to connect as many pieces as desired without losing. In other words---and perhaps to our surprise---P1 can draw the $h=1$ game even if P2 cannot lose, and P2 can draw the $h=1$ game even if P1 cannot lose. 

This approach, where we prove each player can draw the $h=1$ game even if the other cannot lose, is valuable because of its ability to generalize: by showing that one can avoid connecting $3$ even if the other player is unrestricted, it follows that one can avoid connecting $k\geq3$. And, if our Lemma shows that either player can draw an $h=1$ game, we will not only have the result that $h=1$ is a draw for $k\geq3$, but Theorem~\ref{thm-delayed} and Corollary~\ref{cor-delayed} will be complete. We now prove such a Lemma, after introducing some preliminary definitions.\footnote{We note, with some regret, that prior to the peer review process, we were unable to write a more simple proof than the one that follows. Simple non-constructive existence proofs eluded us. So did strategy-mirroring arguments. In many cases, we identified corner cases or holes in simple strategies through simulation. In other cases, we attempted to trick our colleagues into thinking about this problem, and this was moderately successful: the idea to define canonical boards, as we do in this proof, and analyze them as automata came from Ryan D.\ Lewis in 2018, after the rest of the paper was finished. We kept returning to the problem every now and then, but without progress. The proof finally came just after tenure and just before fatherhood, during a magical window when it felt reasonable to reconnect with sticky puzzles. In the end, this ended up being the most fun part of the paper, and the fact that one of our reviewers was subsequently drawn into the problem to come up with an even simpler proof is a delight. May our children one day experience the joy of being stuck on a problem for years and the lovely satisfaction when the solution falls into place.}

The definitions and strategies that follow are introduced without loss of generality from the perspective of P1, who plays as \verb|X|, while P2 plays \verb|O|. After showing how P1 can guarantee a draw as  \verb|X|, we will then show how P2 can do the same.\\

To begin, consider an example game in which each player has played twice, which looks like
$$\verb|X-O---XO--|,$$
with \verb|X| to play. 
While one can certainly view this as a single board of $w=10$ with $n=6$ remaining empty spaces, we'll instead observe that it can be split into a set of three smaller boards as follows 
$$\verb|X-|, \qquad \verb|---X|,\qquad \verb|--|,$$
due to the fact that, because \verb|X| can never connect three across an \verb|O|, each \verb|O| divides the board into two independent boards. From this view, when \verb|X| plays, \verb|X| must choose which board to play in, and which space to play within that board. This approach of choosing a board and then choosing a play will define the eventual optimal strategy.

The multiple-board representation above also means that when \verb|O| plays, we can rewrite the board \verb|O| just played in by splitting it and retaining any resulting boards with playable spaces. Thus, each of \verb|O|'s moves has the potential to turn one board into two. This representation also means that the boundaries of the original board and the \verb|O| pieces are equivalent, leading to the following definition.\\

\noindent{\bf Definition (filler, \verb|F|):} Let \verb|F| denote filler which consists of any number of pieces, with no spaces between them, such that (i) the end pieces are \verb|O|, and (ii) there are no more than two adjacent \verb|X| pieces in the middle. For instance, F could represent \verb|OXXOXO|, or could also represent \verb|OOOO|. Let the boundary of the board also be classified as filler. \\

	The definition of \verb|F| allows us to introduce and classify the types of board that \verb|X| may face. The most simple is $A_n$, consisting of $n$ empty positions bordered by \verb|F|. For instance,  
	$$A_2 = \verb|F--F|, \qquad A_6 = \verb|F------F|\ .$$
Note also that each game starts with board $A_w$. The next most simple is $B_n$, with $n$ empty positions and an \verb|X| on one side, bordered by \verb|F|. For instance,
	$$B_2 = \verb|FX--F|, \qquad B_5 = \verb|F-----XF|\ .$$  
	
These definitions are useful because they allow us to more compactly describe the types of boards that arise during play. For example, 
	$$ \verb|F--F| \px \verb|FX-F| \quad \equiv \quad A_2 \px B_1$$
Note that even if \verb|X| had played in the other available space, we'd still write $A_2 \px B_1$, due to the left-right symmetry. 

With these examples in mind, we now define a set of five canonical boards, labeled $A$, $B$, $C$, $D$, and $E$. For each board type, a subscript denotes the number of empty positions. Each definition implicitly includes its left-to-right reflection because of its irrelevance to strategy, meaning that both \verb|FX--F| and \verb|F--XF| are considered $B_2$, for example.

All canonical boards are bounded by filler on either side, and are defined by the presence/absence and number of \verb|X|s next to the filler, with playable spaces in the middle (Table~\ref{table-boardtypes}). $A$ boards have no \verb|X|s. $B$ boards have an \verb|X| on one side. $C$ boards have an \verb|X| on both sides. $D$ boards have two \verb|X|s on one side. $E$ boards have two \verb|X|s on one side and one \verb|X| on the other side.
\begin{table}[h]
	\renewcommand{\arraystretch}{1.5}
	\begin{tabular}{|c|c|l|l|l|}
		\hline
		board & definition & example & \verb|X| play & \verb|O| products\\
		\hline
		$A$ & \verb|F......F| & $A_4$, \verb|F----F| & $A \px B$ & $A \po A$\\
		$B$ &  \verb|FX.....F|& $B_2$, \verb|FX-F| & $B \px C$ & $B \po A,B$\\
		$C$ & \verb|FX....XF|& $C_3$, \verb|FX---XF| & $C \px E$ & $C \po B$\\
		$D$ & \verb|FXX....F|& $D_2$, \verb|F--XXF| & $D \px E$ & $D \po A, D$\\
		$E$ &  \verb|FXX...XF|& $E_1$, \verb|FXX-XF| & none & $E \po B, D$\\
		\hline
	\end{tabular}
	\cprotect\caption{{\bf Canonical boards and the results of playing each.} The five canonical boards are listed with definitions in which \verb|...| represents an arbitrary number of empty spaces. Examples include variations of both left-right reflections to illustrate the flexibility of the definitions. If playing in a particular board type, the optimal move for \verb|X| is shown. Because we will prove that \verb|X| can lead any game to a draw, we show all possible board that could be produced when \verb|O| plays in each, including the possibility that more than one such product is created were \verb|O| to play in the interior of a board.}
	\label{table-boardtypes}
\end{table}

This set of canonical boards allows us to specify what \verb|X| should do, should they choose to play in each board. When playing an $A_n$, \verb|X| should create a $B_{n-1}$. When playing a $B_n$, \verb|X| should create a $C_{n-1}$. When playing a $D_n$, \verb|X| should create an $E_{n-1}$. We assert that \verb|X| needs no rule for playing $E_n$ because optimal play never involves doing so. These plays are summarized in Table~\ref{table-boardtypes}. We also assert that, while more exotic canonical boards could be defined, our solution need not define them because it avoids them by construction. 

We can also take note of the products resulting from all possible \verb|O| plays. For instance, 
$$\verb|FX--F| \po \verb|FXO-F|\ \ \text{or} \ \ \verb|FX-OF|, $$
meaning that \verb|O| can turn a $B$ into an $A$ or a $B$. Including instances where \verb|O| creates multiple boards, all possible \verb|O| products are also tabulated in Table~\ref{table-boardtypes}.

\begin{figure}[h]
	\includegraphics[width=1.0\linewidth]{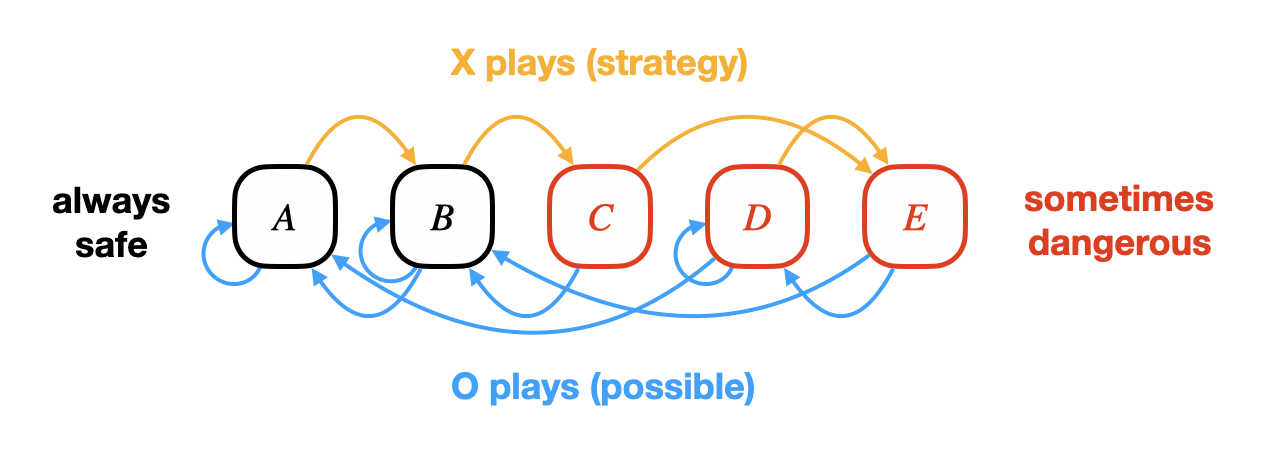}
	\cprotect\caption{{\bf Graphical representation of game flow among canonical board types.} Board types are vertices and plays are directed edges. The specified moves for \verb|X| (orange) and all possible moves for \verb|O| (blue) have broadly different effects: \verb|X|'s plays move strictly to the right and may create boards that are sometimes dangerous (red), while \verb|O|'s tend to (but do not always) move to the left toward safe boards.}
	\label{fig-flow}
\end{figure}

Among the five board types, $A$ and $B$ are always safe for \verb|X| to play in without fear of losing. However, $C$, $D$, and $E$ boards may cause \verb|X| to lose because, if forced to play in them, \verb|X| may have to connect three. In a sense, $A$ and $B$ boards are always safe, and $C$, $D$, and $E$ boards are sometimes dangerous. Now, notice that our specified plays for \verb|X| can never decrease the level of danger, and may increase it. In contrast, many of the possible plays for \verb|O| decrease the level of danger. This flow of plays is represented in Fig.~\ref{fig-flow}, and is meant to create some intuition for how \verb|X| may eventually avoid losing: avoid creating danger when possible, force \verb|O| to turn dangerous boards back into safe ones, and avoid all opportunities for \verb|O| to create new dangers. 

The graphical representation in Fig.~\ref{fig-flow} illustrates an important fact: the set of canonical boards is closed under {\it all possible} \verb|O| plays and {\it all specified} \verb|X| plays. Therefore, as long as we can show that the specified plays lead the game to a draw, we need not consider any other board types. 

We now make one final observation before defining a strategy and proving some results: \verb|X| will lose if forced to play a $C_n$, $D_n$, or $E_n$ board only when $n=1$. It follows that \verb|X| can draw the game if \verb|X| can avoid ever playing in $C_{odd}$, $D_{odd}$, or $E_{odd}$ boards. From this observation, we let ``odd'' and ``even'' subscripts refer to the number of playable spaces in each board, and prove our first Lemma.

\begin{lemma}\label{lemma-odd-pancake}
Suppose that $w$ is odd, $h=1$, the first player \verb|X| loses if they connect $k=3$, and the second player \verb|O| can never lose. The first player \verb|X| can draw the game by playing only in $A_{odd}$ and $B_{odd}$ boards. 
\end{lemma}
\begin{proof}
First, note that if \verb|X| plays only in $A_{odd}$ and $B_{odd}$ boards, then \verb|X| does not lose, and thus---because \verb|O| cannot lose---\verb|X| draws the game. It is therefore sufficient to show only that \verb|X| always has at least one $A_{odd}$ or $B_{odd}$ option.

Second, note that because $w$ is odd and \verb|X| plays first, the total number of spaces available for \verb|X| to play in, summed over all boards, is odd. Any set of positive natural numbers with an odd sum contains at least one odd element, implying that there exists at least one odd board each time \verb|X| must play.

When \verb|X| plays an in-strategy move---that is, when \verb|X| plays only in $A_{odd}$ or $B_{odd}$ boards, as specified in the Lemma we seek to prove---\verb|X| will create either a $B_{even}$ (from an $A_{odd}$) or a $C_{even}$ (from a $B_{odd}$). In short, \verb|X| cannot create a $D$ or $E$ (Fig.~\ref{fig-flow}). Now, observe that, if no $D$ or $E$ boards are available for \verb|O| to play, then it is impossible for \verb|O| to create a $C$, $D$, or $E$ board (again, Fig.~\ref{fig-flow}). Therefore, provided that all of \verb|X|'s moves are in-strategy moves, up until some point in the game, then it is impossible for there to be any $D$ or $E$ boards at all.

Following a similar argument, provided that all of \verb|X|'s moves are in-strategy moves, up until some point in the game, it is impossible for there to be any $C_{odd}$ boards: \verb|O| cannot create them, and \verb|X| creates them only if playing in a $B_{even}$. 

Consequently, as long as \verb|X| plays in-strategy moves, \verb|X| will face only $A$, $B$, and $C_{even}$ boards. However, \verb|X| can never be forced to play $C_{even}$, because of our second introductory observation: there must be at least one odd board every time that \verb|X| plays, and because $C_{even}$ is even, one of the $A$ or $B$ alternatives must be odd. 

Thus, \verb|X| starts the game by playing an $A_{odd}$, and is guaranteed to always have an in-strategy move to play. All in-strategy moves are non-losing. Therefore \verb|X| draws the game.
\end{proof}

Games with even width $w$ require a different strategy, most obviously because the game starts with $A_{even}$ which requires a move that is, by definition, outside the previous proof's strategy of playing only in odd $A$ and $B$ boards. Nevertheless, it remains possible for \verb|X| to draw the game with a simple three-tier prioritization scheme, which we will describe after a few preliminary observations.

\begin{figure*}[t]
	\includegraphics[width=1.0\linewidth]{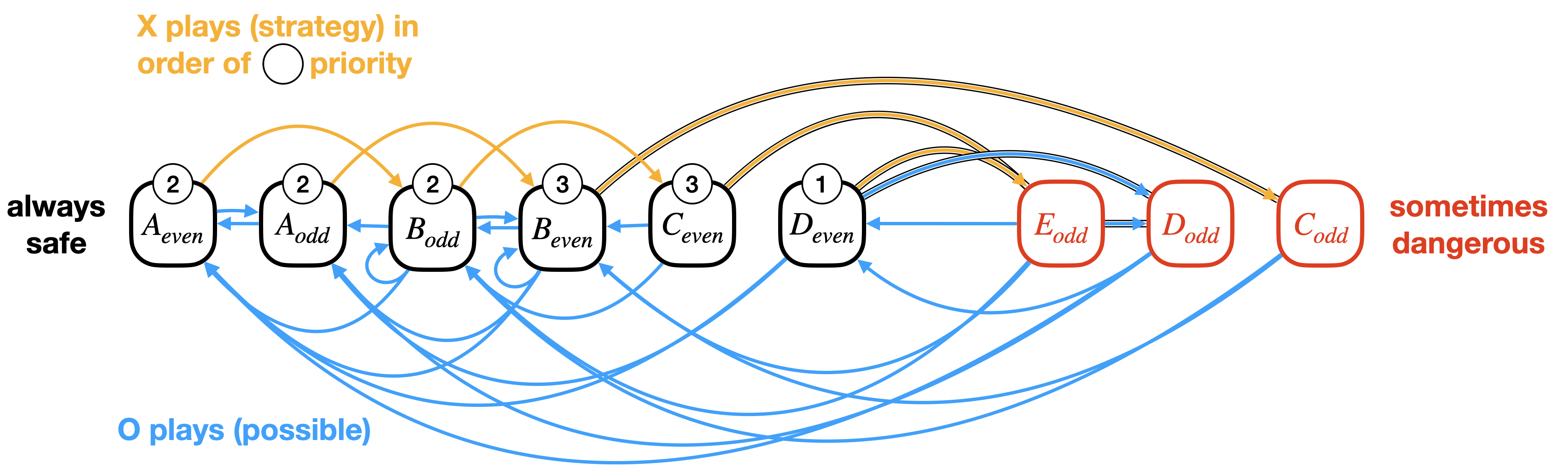}
	\cprotect\caption{{\bf Representation of game flow among even and odd canonical board types.} Board types are vertices and plays are directed edges. The specified moves for \verb|X| (orange) are shown, with annotations for their priority orders among the safe (black) boards. No moves are specified among the dangerous (red) board types because \verb|X| need never play one. All possible moves for \verb|O| (blue) are shown, noting that some moves create two boards and thus follow two edges (e.g. $E_{odd} \po D_{odd} + B_{odd}$). The five moves (three orange and two blue) that can possibly create dangerous boards are bolded with black outlines. The $E_{even}$ board can never be reached by in-strategy play and is not shown.}
	\label{fig-flow-even}
\end{figure*}

Figure~\ref{fig-flow-even} shows a graphical representation of the game flow between canonical board types, this time split into even and odd subtypes. Note that Fig.~\ref{fig-flow-even} contains only nine vertices; $E_{even}$ is missing. This is because only \verb|X| can create an $E_{even}$ from one of the other canonical types, specifically by playing in a $D_{odd}$ next to the filler, e.g.,
$$\verb|FXX---F| \px \verb|FXX--XF|\ .$$
However, as we go on to show in the Lemma that follows, \verb|X| will never be forced to play in a $D_{odd}$. Consequently, because no $E_{even}$ exists at the start of the game, and because neither player will or can create one, the $E_{even}$ board need not clutter our already complicated flow diagram. 

What can we learn about \verb|X|'s strategy from examining Fig.~\ref{fig-flow-even}? First, observe that, although the vast majority of \verb|O|'s possible plays lead to safe boards, there exist two plays that lead to potentially dangerous boards. The first such play occurs when \verb|O| plays into the interior of an $E_{odd}$ board and leaves an odd number of spaces on either side. The side with the \verb|X| becomes a $B_{odd}$ and the side with the \verb|XX| becomes a $D_{odd}$, e.g.
$$\verb|FXX---XF| \po \verb|FXX-O-XF| \quad \equiv\quad \verb|FXX-F|,\quad \verb|F-XF|\ .$$
Importantly, \verb|O| has not increased the number of potentially dangerous boards, but has instead transformed one canonical type of potentially dangerous board into another. 

The second play allowing \verb|O| to produce a dangerous board occurs when \verb|O| plays next to the edge of a $D_{even}$ board to create a $D_{odd}$, e.g.,
$$\verb|FXX----F| \po \verb|FXX---OF| \quad \equiv \quad \verb|FXX---F|\ .$$
In contrast to the first play discussed above, this one allows \verb|O| to increase the number of potentially dangerous boards. This turns out to be very important: if \verb|O| can create more than one dangerous board, \verb|O| can eventually force \verb|X| to play in one of them, and \verb|X| may lose. This means that $D_{even}$ boards are dangerous in the hands of \verb|O|!

The observations above lead to the following counterintuitive strategy: \verb|X| will play $D_{even}$ boards with top priority, just to ensure that \verb|O| cannot---even though doing so means that \verb|X| creates a dangerous $E_{odd}$ board in the process! After prioritizing $D_{even}$ boards, \verb|X| prioritizes safe plays that create no danger ($A_{even}, A_{odd}, B_{odd}$), followed by safe plays that do create danger ($B_{even}, C_{even}$). This simple prioritization, built to disempower \verb|O|, leads to a Lemma for even-width games.

\begin{lemma}\label{lemma-even-pancake}
Suppose that $w$ is even, $h=1$, the first player \verb|X| loses if they connect $k=3$, and the second player \verb|O| can never lose. The first player \verb|X| can draw the game by playing boards in the following order of priority:
\begin{enumerate}[itemsep=0pt,topsep=2pt]
	\item $D_{even}$
	\item $A_{even}$, $A_{odd}$, or $B_{odd}$ 
	\item $B_{even}$ or $C_{even}$
\end{enumerate}
\end{lemma}
\begin{proof}
This proof will show that by playing in the strategic priority, \verb|X| will always have a safe option. In such games, there can be at most one dangerous board, which must always be accompanied by a safe alternative. 

To begin, observe that only \verb|O| can make a $D_{even}$ (Fig.~\ref{fig-flow-even}). While it is possible for \verb|O| to create a lone $D_{even}$, or a single $D_{even}$ and a $B_{even}$, only one $D_{even}$ can be made at a time. Because $D_{even}$ is the one and only Priority 1 board in \verb|X|'s strategy, \verb |X| will always play $D_{even} \px E_{odd}$, which implies that \verb|O| will never see a $D_{even}$. 

The fact that \verb|O| will never see a $D_{even}$ implies that \verb|O| can never transform a safe board into a dangerous board. It can only transform a dangerous $E_{odd}$ board into a dangerous $D_{odd}$ board (Fig.~\ref{fig-flow-even}). In other words, \verb|O| cannot increase the number of dangerous boards. 

If \verb|O| cannot increase the number of dangerous boards directly, then the only way for the number of dangerous boards to increase is if \verb|X| creates one. However, while \verb|X| may be forced to create one dangerous board, it will never create a second: dangerous boards are odd, \verb|X| must face an even total number of playable spaces, so \verb|X| must also have a safe {\it odd} board among their options. In other words, when \verb|X| faces a single dangerous board, they always have a safe (and odd) alternative. Because that safe alternative is guaranteed to be odd, it cannot be a $D_{even}$, and thus \verb|X| will never create a second dangerous board.

We have shown that the strategic priority limits \verb|O| from creating dangerous boards directly, and they can force \verb|X| to create at most one dangerous board. Yet that single dangerous board must be accompanied by a safe alternative. Therefore \verb|X| draws the game.
\end{proof}

These results now lead us to the following Theorem, which extends the results above for $k=3$ to arbitrary $k$, and from P1 to P2.

\begin{theorem}
	Suppose that $h=1$, $w$ is arbitrary but finite, and either player loses if they connect $k\geq3$. Perfect play allows either player to draw the game.
\end{theorem}
\begin{proof}
	Lemmata~\ref{lemma-odd-pancake} and \ref{lemma-even-pancake} show that P1 can avoid connecting $k=3$. Importantly, these proofs allow P2 to play arbitrarily, meaning that P2 could have connected 3 or more, yet P1 can nevertheless bring the game to a draw. If P2 is entirely unrestricted in their options, and yet P1 can avoid connecting 3, then it follows that P1 can avoid connecting $k\geq3$. Thus, P1 can draw the game for arbitrary finite width $w$ and arbitrary $k\geq3$.
	
	Now consider the game from P2's perspective. In the game's opening move, P1 will either play against a boundary or in the empty board's interior. On the one hand, if P1 chooses to open the game by playing against a boundary, then P2 faces a safe $A_{w-1}$ board. Letting P2 play \verb|X| (simply to use the previously established canonical board types), the previous Lemmata apply, and P2 can draw the game.
	
	On the other hand, if P1 chooses to open the game by playing in the interior, then P2 faces two boards, $A_u$ and $A_v$ with $u+v = w-1$. If $w$ is even then $w-1$ is odd, meaning that either $u$ or $v$ is odd, and Lemma~\ref{lemma-odd-pancake} applies. If $w$ is odd, then $w-1$ is even, meaning that either $u$ and $v$ are both odd, or $u$ and $v$ are both even. Because these two even or two odd boards are both type $A$, the boards are not dangerous, nor can play therein by P2 create a dangerous board, and thus Lemma~\ref{lemma-even-pancake} applies, and P2 can draw this type of game as well.
	
	Because both players can avoid losing, \mck with $h=1$ and $k\geq3$ is a draw.
\end{proof}

%

\end{document}